\documentclass [a4paper,twoside,12pt]{article}

\usepackage{amsmath,amsfonts,amssymb}
\usepackage{vmargin,graphicx,theorem}
\usepackage{enumerate}

\setpapersize[portrait]{A4}
\setmarginsrb{1cm}{2.2cm}{1cm}{2.5cm}{0cm}{0.1cm}{0.5cm}{2cm}

\newcommand{\rr}{\ensuremath{\mathbb{R}}}

\newtheorem{ethm}{}
\newtheorem{theo}[ethm]{Theorem}

\newtheorem{erem}[ethm]{Remark}

\newcommand{\proofend}{~$\rhd$}
\newcommand{\proofbegin}{~$\lhd$}

\newcommand{\p}[4]{{#3}\!\left#1{#4}\right#2}

\newcommand{\PAR}[1]{\ensuremath{{\left(#1\right)}}} 
\newcommand{\SBRA}[1]{\ensuremath{{\left[#1\right]}}} 

\renewcommand{\phi}{\varphi}

\newcommand{\varf}[1]{\mathbf{Var}_{#1}}
\newcommand{\entf}[1]{\mathbf{Ent}_{#1}}

\newcommand{\ent}[2]{\p(){\entf{#1}}{#2}}

\newcommand{\var}[2]{\p(){\varf{#1}}{#2}}

\newcommand{\PT}[1]{\mathbf{P}_{\!#1}}
\newcommand{\Pt}[1][t]{\ensuremath{\mathbf{P_{\!#1}}}}

\begin{document}

\title{\sl Phi-entropy inequalities and Fokker-Planck equations}
\author{ Fran\c cois Bolley and Ivan Gentil }

\date{November, 2010 }
\maketitle\thispagestyle{empty}

\bigskip

\abstract{We present new $\Phi$-entropy inequalities for diffusion semigroups under the curvature-dimension criterion. They include the isoperimetric function of the Gaussian measure. Applications to the long time behaviour of solutions to Fokker-Planck equations are given.}

\bigskip

\bigskip



\bigskip

\section{Introduction}

We consider a Markov semigroup $(\PT{t})_{t \geq 0}$ on $\rr^n$, acting on functions on $\rr^n$ by
$\PT{t} f(x) =  \int_{\rr^n} f(y) \, p_t(x, dy)$
for $x$ in $\rr^n.$ The kernels $p_t(x, dy)$ are probability measures on $\rr^n$ for all $x$ and $t \geq 0$, called transition kernels. We assume that the Markov infinitesimal generator $L = \displaystyle \frac{\partial  \PT{t}}{\partial t} \Big\vert_{t=0^+}$ is given~by
$$
L f (x)= \sum_{i,j=1}^n D_{ij} (x) \frac{\partial^2 f}{\partial x_i \partial x_j} (x) - \sum_{i=1}^n a_i (x)  \frac{\partial f}{\partial x_i} (x)
$$
where  $D(x) = (D_{ij}(x))_{1 \leq i,j \leq n}$ is a symmetric $n \times n$ matrix, nonnegative in the sense of quadratic forms on $\rr^n$ and with smooth coefficients, and where the $a_i, 1 \leq i \leq n,$ are smooth. Such a semigroup or generator is called a {\it diffusion}, and we refer to Refs.~\cite{bakrystflour},  \cite {bakrytata},~\cite{ledouxmarkov}  for backgrounds on them. 

If $\mu$ is a Borel probability measure on $\rr^n$ and $f$ a $\mu$-integrable map on $\rr^n$ we let $\mu(f) = \int_{\rr^n} f(x) \, \mu(dx).$
If, moreover, $\Phi$ is a convex map on an interval $I$ of $\rr$ and $f$ an $I$-valued map with $f$ and $\Phi(f)$ $\mu$-integrable,  we let 
$$
\entf{\mu}^\Phi\PAR{f} = \mu(\Phi(f)) - \Phi(\mu(f)) 
$$
be the $\Phi$-entropy of $f$ under $\mu$ (see Ref.~\cite{chafai04} for instance).
Two fundamental examples are  $\Phi(x) = x^2$ on $\rr$, for which $\entf{\mu}^\Phi\PAR{f}$ is the variance of $f,$ and $\Phi(x) = x \ln x $ on $]0, +\infty[$, for which $\entf{\mu}^\Phi\PAR{f}$ is the Boltzmann entropy of $f$.
 By Jensen's inequality, $\entf{\mu}^\Phi (f)$ is always nonnegative and, if $\Phi$ is strictly convex, it is positive unless $f$ is a constant, equal to $\mu(f).$ 
 The semigroup $(\Pt)_{t \geq 0}$ is said {\it $\mu$-ergodic} if $\Pt f$ tends to $\mu(f)$ as $t$ tends to infinity in $L^2(\mu)$, for all~$f$. 
 
 In Section~\ref{sectone} we shall derive bounds on  $\entf{\mu}^\Phi\PAR{f}$ and  $\entf{\PT{t}}^\Phi\PAR{f} (x)$ which will measure the convergence of $\Pt f$ to $\mu(f)$ in the ergodic setting. This is motivated by the study of the long time behaviour of solutions to Fokker-Planck equations, which will be discussed in Section~\ref{secttwo}. 
 
Some results of this note with their proofs are detailled in Ref.~\cite{bolley-gentil}.

\section{Phi-entropy inequalities}\label{sectone}
Bounds on $\entf{\PT{t}}^\Phi\PAR{f}$ and assumptions on $L$ will be given in terms of the {\it carr\'e du champ}  and  $\Gamma_2$ operators associated to $L$, defined  by
 $$
 \Gamma(f,g)= \frac{1}{2} \Big( L(fg) - f \, Lg - g \, Lf \Big), \,\,  \Gamma_2(f) = \frac{1}{2} \Big( L \Gamma(f) - 2\Gamma (f, Lf) \Big).
 $$
 
If $\rho$ is a real number, we say that the semigroup $(\PT{t})_{t \geq 0}$ satisfies the {\it $CD(\rho,\infty)$ curvature-dimension} (or Bakry-\'Emery) {\it criterion} (see Ref.~\cite{bakryemery})  if
$$
\Gamma_2(f)\geq \rho \, \Gamma(f)
$$
for all functions $f$, where $\Gamma (f) = \Gamma (f,f).$

The carr\'e du champ is explicitely given by 
$$\Gamma (f,g) (x) =  < \nabla f (x), D(x) \,  \nabla g(x) >.$$
Expressing $\Gamma_2$ is more complex in the general case but, for instance, if $D$ is constant, then  $L$ satisfies the $CD(\rho,\infty)$  criterion if and only if
\begin{equation}\label{cdrhocst}
\frac{1}{2} \big( J a (x) D + (J a (x) D )^* \big) \geq \rho \, D
\end{equation}
for all $x$, as quadratic forms on $\rr^n$, where $\displaystyle J a$ is the Jacobian matrix of $a$ and $M^*$ denotes the transposed matrix of a matrix $M$  (see Ref.~\cite{arnoldcarlenju08,amtucpde01}) .

%
 
 Poincar\'e and logarithmic Sobolev inequalities for the semigroup $(\PT{t})_{t \geq 0}$ are known to be implied by the ${CD}(\rho,\infty)$ criterion. 
  More generally, and following Ref.~\cite{bakrytata,bakryemery,chafai04}, let $\rho>0$ and $\Phi$ be a $C^4$ strictly convex function on an  interval $I$ of $\rr$ such that $-1/\Phi''$ is convex.
 If $(\Pt)_{t \geq 0}$ is $\mu$-ergodic and satisfies the $CD(\rho, \infty)$ criterion,  then $\mu$ satisfies the $\Phi$-entropy inequality 
\begin{equation}\label{PHII}
 \entf{\mu}^\Phi (f)   \leq \frac{1}{2 \rho} \, \mu (\Phi''(f) \, \Gamma (f))
\end{equation}
for all $I$-valued functions $f$.

The main instances of such $\Phi$'s are the maps $x \mapsto x^2$ on $\rr$ and $x \mapsto x \ln x$  on $]0, +\infty[$ or more generally, for $1 \leq p \leq 2$ 
\begin{equation}
\label{eq-phip}
\Phi_p(x)= \left\{
\begin{array}{cl}
 \frac{x^p- x}{p(p-1)}, \quad x > 0 & \text{ if }p\in ]1,2]\\
 x \ln x, \quad x>0 & \text{ if }p=1.
\end{array}
\right.
\end{equation}
For this $\Phi_p$ with $p$ in $]1,2]$ the $\Phi$-entropy inequality \eqref{PHII}  becomes
\begin{equation}\label{beckner}
\frac{\mu(g^2) - \mu(g^{2/p})^p}{p-1} \leq \frac{2}{ p  \rho} \, \mu (\Gamma (g))
 \end{equation}
for all positive functions $g$. For given $g$ the map $p \mapsto \frac{\mu(g^2) - \mu(g^{2/p})^p}{p-1}$ is nonincreasing with respect to $p >0, p\neq1$. Moreover its limit for $p \to 1$ is $\ent{\mu}{g^2}$, so that
 the so-called {\it Beckner inequalities}~\eqref{beckner} for $p$ in $]1,2]$ give a natural monotone interpolation between the weaker Poincar\'e inequality (for $p=2$), and the stronger logarithmic Sobolev inequality (for $p \to 1$).

\subsubsection*{Long time behaviour of the semigroup}
The $\Phi$-entropy inequalities provide estimates on the long time behaviour of the associated diffusion semigroups. 
Indeed, let $(\Pt )_{t\geq 0}$ be such a semigroup, ergodic for the measure $\mu$.  If $\Phi$ is a $\mathcal C^2$ function on an interval $I$,~then 
\begin{equation}\label{HI}
\frac{d}{dt}\entf{\mu}^\Phi\PAR{\Pt f}=-\mu\PAR{\Phi''(\Pt f) \, \Gamma(\Pt f)}
\end{equation}
for all $t \geq 0$ and all $I$-valued functions $f.$ As a consequence, if $C$ is a positive number, then the semigroup converges in $\Phi$-entropy with exponential rate: 
\begin{equation}\label{cv1}
\entf{\mu}^\Phi\PAR{\Pt f}\leq e^{-\frac{t}{C}}\entf{\mu}^\Phi\PAR{ f}
\end{equation}
for all $t \geq 0$ and all $I$-valued functions $f$, if and only if the measure $\mu$ satisfies the $\Phi$-entropy inequality  for all $I$-valued functions $f$,
\begin{equation}
\label{eq-phisob}
\entf{\mu}^\Phi\PAR{f}\leq C\mu \PAR{\Phi''( f) \, \Gamma(f)}.
\end{equation}

\subsection{Refined $\Phi$-entropy inequalities}

We now give and study improvements of~\eqref{PHII} for the  $\Phi_p$ maps given by~\eqref{eq-phip}:
\begin{theo}[{[\cite{bolley-gentil}]}]
\label{thm-main}
Let $\rho \in \rr$ and  $p \in ]1,2[$. Then the following assertions are equivalent, with $\PAR{1 - e^{-2\rho t}}/{\rho}$ and $\PAR{e^{2\rho t} -1}/{\rho}$ replaced by $2 t$ if $\rho =0$:
\begin{enumerate}
\item[(i)]
 the semigroup $(\PT{t})_{t \geq 0}$ satisfies the ${CD}(\rho,\infty)$ criterion;
\item [(ii)]
 $(\PT{t})_{t \geq 0}$ satisfies the refined local $\Phi_p$-entropy inequality
$$
\frac{1}{(p-1)^2}\SBRA{{\Pt(f^p)}-{\Pt(f)^p}\PAR{\frac{\Pt(f^p)}{\Pt(f)^p}}^{\! \! \frac{2}{p}-1}}\leq \frac{1-e^{-2\rho t}}{\rho}\Pt\PAR{f^{p-2} \, {\Gamma(f)}}
$$
 for all positive $t$ and all positive functions $f$;
 \item[(iii)]
$(\PT{t})_{t \geq 0}$ satisfies the reverse refined local $\Phi_p$-entropy inequality  
$$
\frac{1}{(p-1)^2}\SBRA{{\Pt(f^p)}-{\Pt(f)^p}\PAR{\frac{\Pt(f^p)}{\Pt(f)^p}}^{\! \! \frac{2}{p}\!-\!1}}
\geq
 \frac{e^{2\rho t}\! \!-\!1}{\rho}\PAR{\frac{(\Pt f)^p}{\Pt\PAR{f^p}}}^{\! \! \frac{2}{p}\!-\!1}\!\!\!\PAR{\Pt f}^{p-\!2} \, {\Gamma(\Pt f)}
$$
 for all positive $t$ and all positive functions $f$.
  \end{enumerate}
   If, moreover, $\rho >0$ and the measure $\mu$ is ergodic for the semigroup $(\Pt)_{t \geq 0}$, then $\mu$ satisfies the refined $\Phi_p$-entropy inequality
\begin{equation}
\label{eq-ad2}
\frac{p^2}{(p-1)^2}\SBRA{{\mu(g^2)}-{\mu(g^{2/p})^p}\PAR{\frac{\mu(g^2)}{\mu(g^{2/p})^p}}^{\frac{2}{p}-1}}\leq \frac{4}{\rho}\mu\PAR{\Gamma(g)}
\end{equation}
for all positive maps $g.$

\end{theo}

The bound~\eqref{eq-ad2} has been obtained  in Ref.~\cite{arnolddolbeault05}  for the generator $L$ defined by $Lf = \textrm{div} (D \nabla f) - < \! \!D \nabla V , \nabla f \! \!>$ with $D(x)$ a scalar matrix and for the ergodic measure $\mu = e^{-V}$, and under the corresponding $CD(\rho, \infty)$ criterion. 

It improves on the Beckner inequality~\eqref{beckner} since 
\begin{equation}\label{comp-beckner+}
\frac{\mu(g^2) - \mu(g^{2/p})^p}{p-1} \leq \frac{p}{2 (p-1)^2} \Big[ \mu(g^2) -\mu(g^{2/p})^p \, \Big(\frac{\mu(g^2)}{\mu(g^{2/p})^p}\Big)^{\frac{2}{p}-1}  \Big].
\end{equation}
We have noticed that for all $g$ the map 
$
p \mapsto \frac{\mu(g^2) - \mu(g^{2/p})^p}{p-1} 
$
is continuous and nonincreasing on $]0, +\infty[$, with values $\ent{\mu}{g^2}$ at~$p=1$ and $\var{\mu}{g}$ at $p=2$. Similarly, for the larger  functional introduced in \eqref{comp-beckner+},~the~map 
$$
p \mapsto \frac{p}{2 (p-1)^2} \Big[ \mu(g^2) -\mu(g^{2/p})^p \, \Big(\frac{\mu(g^2)}{\mu(g^{2/p})^p}\Big)^{\frac{2}{p}-1}  \Big]
$$ 
is nonincreasing on $]1,+\infty[$ (see~\cite[Prop. 11]{bolley-gentil}). Moreover its value is $\var{\mu}{g}$  at $p=2$ and it tends to $ \ent{\mu}{g^2}$ as $p \to 1,$ hence providing a new monotone interpolation between Poincar\'e and logarithmic Sobolev inequalities.

\medskip

The pointwise $CD(\rho, \infty)$ criterion can be replaced by the integral criterion 
$$
\mu\PAR{g^{\frac{2-p}{p-1}}\Gamma_2(g)}\geq  \rho \, \mu\PAR{g^{\frac{2-p}{p-1}}\Gamma(g)}
$$
for all positive functions $g,$  and one can still get the refined $\Phi_p$-entropy inequality~\eqref{PHII}, even in the case of non-reversible semigroups (see~\cite[Prop.~14]{bolley-gentil}).

 \begin{erem}
 For $\rho=0$, and following Ref.~\cite{arnolddolbeault05}, the convergence of $\Pt f$ towards $\mu(f)$ can be measured on $H(t) = \entf{\mu}^\Phi\PAR{\Pt f}$ as
$$
\vert H'(t)\vert \leq \frac{\vert H'(0) \vert}{1+ \alpha t}, \quad t \geq 0
 $$
 where $\alpha = \frac{2-p}{p} \vert H'(0) \vert / H(0).$ This  illustrates the improvement offered by~\eqref{eq-ad2} instead of~\eqref{beckner}, which does not give here any convergence rate.
 \end{erem}

\subsection{The case of the Gaussian isoperimetry function}

Let $F$ be the distribution function of the one-dimensional standard Gaussian measure. 
The map $\mathcal U=F' \circ F^{-1}$, which is the isoperimetry function of the Gaussian distribution, satisfies $\mathcal U''=-1/\mathcal U$ on the set $[0,1]$, so that the map $\Phi=-\mathcal U$ is convex with $-1/\Phi''$ also convex on $[0,1].$
 
\begin{theo}\label{theophisoperimetry}
Let $\rho$ be a real number. Then the following three assertions are equivalent, with $\PAR{1 - e^{-2\rho t}}/{ \rho}$ and $\PAR{e^{2\rho t} -1}/{ \rho}$ replaced by $2t$ if $\rho =0$:
\begin{enumerate}
\item [(i)]
 the semigroup $(\PT{t})_{t \geq 0}$ satisfies the ${CD}(\rho,\infty)$ criterion;
\item [(ii)]
 the semigroup $(\PT{t})_{t \geq 0}$ satisfies the local $\Phi$-entropy inequality
 \begin{equation*}
 \entf{\Pt}^\Phi (f) \leq \frac{1}{\Phi''(\Pt f)}\log\PAR{1+\frac{1 - e^{-2\rho t}}{2 \, \rho}\Phi''(\Pt f) \,  \Pt( \Phi''(f) \Gamma (f))}
 \end{equation*}
   for all positive $t$ and all $[0,1]$-valued functions $f$; 
 \item[(iii)]
the semigroup $(\PT{t})_{t \geq 0}$ satisfies the reverse local $\Phi$-entropy inequality
  \begin{equation*}
  \entf{\Pt}^\Phi (f) \geq \frac{1}{\Phi''(\Pt f)}\log\PAR{1+ \frac{e^{2\rho t} -1}{2 \, \rho}\Phi''(\Pt f)^2 \Gamma (\Pt f))}
  \end{equation*}
    for all positive $t$ and all $[0,1]$-valued functions $f$.
 \end{enumerate}
 If, moreover, $\rho>0$ and the measure $\mu$ is ergodic for the semigroup $(\Pt)_{t \geq 0}$, then $\mu$ satisfies the $\Phi$-entropy inequality for all $[0,1]$-valued functions $f$:
\begin{equation*}
 \entf{\mu}^\Phi (f)   \leq  \frac{1}{\Phi''(\mu( f))}\log\PAR{1+\frac{\Phi''(\mu( f))}{2 \rho} \mu( \Phi''( f)\Gamma (f))}. 
\end{equation*}
\end{theo}

The proof is based on~\cite[Lemma~4]{bolley-gentil}. For $\Phi = - \mathcal U$ it  improves on  the general $\Phi$-entropy inequality~\eqref{PHII} since  $\log(1+x)\leq x$.  Links with the isoperimetric bounds of Ref.~\cite{bakryledoux} for instance will be addressed elsewhere.

\section{Long time behaviour for Fokker-Planck equations}\label{secttwo}
Let us consider the linear Fokker-Planck equation
\begin{equation}
\label{eq-fp1}
\frac{\partial u_t}{\partial t}  = \textrm{div} \SBRA{  D(x)( \nabla u_t + u_t (\nabla V(x)+F(x)))}, \quad t \geq 0, \, x \in \rr^n
\end{equation}
where $D(x)$ is a positive symmetric $n \times n$ matrix and $F$ satisfies 
\begin{equation}\label{rev}
\textrm{div} \PAR{e^{-V} D F}=0.
\end{equation}
It is one of the purposes of Refs.~\cite{arnoldcarlenju08} and~\cite{amtucpde01} to rigorously study the asymptotic behaviour of solutions to~\eqref{eq-fp1}-\eqref{rev}.  Let us formally rephrase the argument.  

Assume that the Markov diffusion generator $L$ defined by 
 \begin{equation}\label{fparnold2}
Lf  =  \textrm{div} (D \nabla  f)  -  < D ( \nabla V-F) , \nabla f>
\end{equation}
 satisfies the $CD(\rho,\infty)$ criterion with $\rho>0$, 
that is~\eqref{cdrhocst} if $D$ is constant, etc.

Then the semigroup $(\Pt)_{t \geq 0}$ associated to $L$ is $\mu$-ergodic with $d\mu=e^{-V}/Zdx$ where $Z$ is a normalization constant. Moreover, a $\Phi$-entropy inequality~\eqref{eq-phisob} holds with  $C=1/(2\rho)$ by~\eqref{PHII}, so that the semigroup converges to $\mu$ according to~\eqref{cv1}. 
However, under~\eqref{rev}, the solution to~\eqref{eq-fp1} for the initial datum $u_0$ is given by $u_t=e^{-V}\,\Pt (e^V u_0)$. Then we can deduce the convergence of the solution $u_t$ towards the stationary state $e^{-V}$ (up to a constant) from the convergence estimate~\eqref{cv1} for the semigroup, in the form
\begin{equation}\label{cvedp}
\entf{\mu}^\Phi\PAR{\frac{u_t}{e^{-V}}}\leq e^{-2\rho{t}}\entf{\mu}^\Phi\PAR{ \frac{u_0}{e^{-V}}}, \quad t \geq 0.
\end{equation}

\bigskip

In fact such a result holds for the general Fokker-Planck equation
\begin{equation}
\label{eq-fp2}
\frac{\partial u_t}{\partial t}  = \textrm{div} \SBRA{  D(x) (\nabla u_t + u_t a(x) )}, \quad t \geq 0, \, x \in \rr^n
\end{equation}
where again $D(x)$ is a positive symmetric $n \times n$ matrix and $a(x) \in \rr^n.$ Its generator is the dual (for the Lebesgue measure) of the generator 
\begin{equation}\label{exam}
Lf = \textrm{div} (D \nabla  f)  -  < D a , \nabla f>. 
\end{equation}
Assume that the semigroup associated to $L$ is ergodic and that its invariant probability measure $\mu$ satisfies a $\Phi$-entropy inequality~\eqref{eq-phisob} with a constant $C\geq0$: this holds for instance if $L$ satisfies the $CD(1/(2C),\infty)$ criterion.

In this setting when $a (x)$ is not a gradient, the invariant measure $\mu$ is not explicit. Moreover the relation  $u_t=e^{-V}\,\Pt (e^V u_0)$ between the solution of~\eqref{eq-fp2} and the semigroup associated to $L$ does not hold, so that the asymptotic behaviour~\eqref{cvedp} for solutions to~\eqref{eq-fp2}  can not be proved by using~\eqref{cv1}. 
However, this relation can be replaced by the following argument, for which the ergodic measure is only assumed to have a positive density $u_{\infty}$ with respect to the Lebesgue measure.

Let $u$ be a solution of~\eqref{eq-fp2} with initial datum $u_0$. Then, by~\cite[Lemma~7]{bolley-gentil},
$$
\frac{d}{dt}\entf{\mu}^\Phi\big( \frac{u_t}{u_\infty}  \big) = \int  \Phi'  \big ( \frac{u_t}{u_\infty}  \big)L^* u_t  dx
=\int  L \Big[  \Phi'  \big(\frac{u_t}{u_\infty} \big) \Big]  \frac{u_t}{u_\infty} d\mu
=  -  \int  \Phi'' \big(\frac{u_t}{u_\infty} \big)\Gamma \big( \frac{u_t}{u_\infty} \big)d\mu.
$$
 Then a $\Phi$-Entropy inequality~\eqref{eq-phisob} for $\mu$ implies the exponential convergence:
\begin{theo}
With the above notation, assume that a $\Phi$-entropy inequality~\eqref{eq-phisob} holds for   $\mu$ and with a constant $C$.  Then all solutions $u = (u_t)_{t \geq 0}$   to the Fokker-Planck equation \eqref{eq-fp2} converge to $u_{\infty}$ in $\Phi$-entropy, with
$$
\entf{\mu}^\Phi\PAR{\frac{u_t}{u_\infty}}\leq e^{-t/C}\entf{\mu}^\Phi\PAR{ \frac{u_0}{u_\infty}}, \quad t \geq 0.
$$
\end{theo}



\noindent
{\bf Acknowledgment:}
This work was presented during the 7th ISAAC conference held in Imperial College, London in July 2009. It is a pleasure to thank the organizers for giving us this opportunity.

\end{document}